# Quantum Reversal of Dialectics: Inversion of Qualities and Occupation Numbers in Statistical Physics

Elmira. A. Isayeva

Institute of Physics of the Ministry of Science and Education of the Republic of Azerbaijan, Baku, Az-1143, Javid av.33, e-mail: el_max63@yahoo.com

## Abstract

This work proposes an original conceptual approach to the derivation and interpretation of quantum statistics — the method of inversion of qualities. Based on philosophical ideas of dialectics (the reversed Hegelian triad of the transition of quality into quantity), W. Ostwald's energetics, and concepts of modern superstring theory, the fundamental dualism of quantum particles and states is reconsidered. Instead of the classical distribution of distinguishable objects over energy levels, we identify the quantum state itself as the primary Quality, while its occupation number (the number of particles) is treated as Quantity.

It is shown that applying combinations without repetition to the simplest two-level system directly gives rise to fermionic antisymmetry and the Pauli exclusion principle, which underlies the structure of the periodic table. In turn, applying combinations with repetition to equidistant systems (the quantum harmonic oscillator) provides a direct and elegant derivation of the bosonic denominator in M. Planck's formula for black-body radiation. The limits of applicability of the method are outlined by contrasting phonon modes in the Einstein crystal model with the nonequidistant spectrum of the hydrogen atom. It is demonstrated that the method of inversion of qualities is conceptually identical to Dirac's second-quantization formalism in modern quantum field theory.

## Introduction: From Hegel and Ostwald to String Theory

Classical German philosophy, in the person of Hegel, formulated one of the principal laws of development: the law of the transition of quantity into quality [1]. In the macroscopic world this principle appears unshakable: we take discrete, tangible objects (quantity) and, by accumulating them, obtain new properties of a system (quality) — thus the chaotic motion of trillions of molecules gives rise to thermodynamic temperature, while heating water to 100 degrees turns it into steam. Everyday experience has deeply rooted in us the illusion that "things" are primary, while their states are secondary.

At the end of the nineteenth century, the German physical chemist Wilhelm Ostwald [3] made a radical attempt to reconsider this materialistic view by developing the philosophy of energetics. Ostwald argued that matter is a fiction, an illusion of our perception, and that in reality there are no "atoms" as indivisible solid building blocks. "Everything we know about the world is only manifestations of energy," he proclaimed. From the standpoint of energetics, what we call matter is merely a stable node of intersecting energy flows. The physics of the time rejected Ostwald's ideas; however, the quantum revolution of the twentieth century showed that his intuition was remarkably close to the truth.

Today, at the frontier of theoretical physics, superstring theory [12] carries this idea to its ultimate conclusion. In its conception there are no elementary particles as pointlike objects. All the diversity of the Universe — electrons, quarks, photons — consists only of different modes of vibration of unified fundamental quantum strings. Just as the same violin string, under different tension and at different vibration frequencies, produces the notes "C," "E," or "G," a superstring through its vibration produces mass, charge, and particle type.

Here the circle closes, and Hegel's classical dialectics undergoes a quantum reversal. At the microscopic level of the Universe, Quality passes into Quantity [1, 10].

The fundamental and primary element of reality becomes not the "particle" (object), but its quantum state — Quality (a definite mode of vibration of a field or string). The "particle" itself arises only as Quantity — as the degree of excitation (occupation number) of this quality.

In this work we propose a method of inversion of qualities in quantum combinatorics. Combinatorial methods have long been used to derive quantum statistics. In the classical approach, the ways of distributing particles among quantum states are counted. In the present work, attention is shifted to the states themselves and their occupation, which makes it possible to reinterpret the combinatorial procedure. We interchange the familiar classical entities: we abandon the counting of "objects distributed over levels" and move to counting "levels possessing a definite quantum quality." This approach not only simplifies the mathematical apparatus, but also reveals a common logical root from which the fundamentally different worlds of fermions and bosons emerge. Let us recall that Dirac was the first to begin counting not particles, but quantum states themselves. This is why his occupation numbers describe our inversion of quality and quantity so naturally and are conceptually consistent with Dirac's second-quantization formalism. But it is not only Dirac. If we recall Fock spaces, creation and annihilation operators, and the way physicists speak about orbitals, field modes, and level occupation numbers, we see that modern quantum statistics indeed works with quantum states and their occupation numbers. In this sense, attention to states is not new. Understanding the identity of microscopic objects led physicists away from Gibbs's permutations with repetition and Boltzmann's combinations without repetition and toward Bose statistics, in which Planck's combinations with repetition are used.

A very simple explanation of why this happens is as follows. Because the microscopic world is invisible to us, its objects are identical to one another. Therefore, their permutations and arrangements have no physical meaning. In a measurement, a detector registers not the individuality of objects, but their action, the energy transferred by them, first represented by Planck as equal portions of energy — energy quanta. However, in the statistics of this action, Boltzmann, mathematically using combinations without repetition and counting microstates, makes his choice among distinguishable objects. By making a "choice," the object is already distinguished and marked — Boltzmann, as it were, points a "finger" at the object: "this one" or "not this one." But for Planck all objects are identical. If Boltzmann answers the question, "Which object possesses this energy?", Planck answers the question, "How many quanta of energy are present in the object, in the system?" This leads to combinations with repetition. Planck points his "finger" not at the object (he does not "see" it in the microscopic world), nor at the quantity, but at the action of the object on the detector, at the quality, representing this action by Planck's constant. If this action occurs frequently, its energy is large — it is one kind of action, one kind of energy quantum. If it occurs rarely, it is another energy quantum. There are different frequencies, and therefore different energy quanta acting on the detector, and their counting begins: 1 such quantum, 2 such quanta, and so on, R such quanta. The problem then becomes one of distributing these R energy quanta among N identical objects. This gives rise to unusual unordered sets of R elements — energy quanta (which have the meaning of quantities rather than qualities), chosen from N types of elements — the microscopic object itself, behind which stand Planck's resonators (which have the meaning of qualities rather than quantities), leading Planck to combinations with repetition. $\varepsilon(1\varepsilon)(2\varepsilon)(P\varepsilon),\ that\ is\ \ A_P^N = \frac{(N+P-1)!}{(N-1)!P!}$

To explain the proposed inversion, let us consider a simple analogy. In the classical description, a room may contain different, individually distinguishable objects: for example, five chairs and two tables. A chair and a table represent different qualities, while individual objects can be directly counted: 5 chairs, 2 tables. At the same time, repetition of the same quality is possible: all seven objects may be chairs.

In the quantum description, the proposed inversion changes the very principle of counting. A quantum object is regarded as an integral state within which different modes — different kinds of quanta — exist. Individual identical particles are not distinguished as individual objects. The states (modes) are distinguished, while the number of their realizations is specified by the occupation number.

If a quantum of definite energy and frequency corresponds to one mode, then one action registered by the detector corresponds to one quantum of that mode. Repetition of the same action means an increase in the number of quanta: 1, 2, 3, … while the same mode is preserved. Thus, the detector registers the action of a quantum of definite energy, while the number of registered quanta characterizes the quantitative aspect of this action.

In this sense, an inversion of the relation between quality and quantity takes place: states (modes) become qualitatively distinguishable, whereas identical particles act as their quantitative realizations. Therefore, the state of a quantum system is naturally specified by a set of occupation numbers: $(n_1, n_2, \ldots, n_r)$, where $n_i$ is the number of quanta in the i-th mode.

**Geometric Interpretation of the Space of Modes**

In the proposed approach, the quantum object is regarded as a space of possible states in which each distinguishable mode corresponds to a separate degree of freedom. For illustration, two modes can be associated with two coordinate axes: the first axis characterizes the number of quanta in the first mode, and the second the number of quanta in the second mode. In general, each mode $\nu$ corresponds to a separate direction in the space of occupation numbers.

The state of the system is specified by the set of occupation numbers $(n_1, n_2, \ldots, n_r)$,

where $n_\nu$ determines the number of quanta occupying the $\nu$-th mode. Thus, occupation numbers act as coordinates of the state, and the collection of all admissible sets $(n_1, n_2, \ldots, n_r)$ forms the space of possible quantitative states of the system.

In the two-dimensional case, with a fixed total number of quanta N

$n_1 + n_2 = N$

the admissible states form the sequence of points

$(0, N), (1, N-1), \ldots, (N, 0)$.

Each point corresponds to a definite distribution of the number of quanta between the two modes, whereas the quanta belonging to one mode have no individual distinctions. Consequently, what is physically significant is not the individual "address" of a particle, but its belonging to a particular mode, that is, its position in the space of modes.

In this sense, one may speak of a quantum's "address": its position in state space is determined by the mode, while the quantitative characteristic of that mode is given by its occupation number. Such an "address" is not a spatial coordinate of a particle in ordinary three-dimensional space; it is a coordinate in the space of quantum states.

It is precisely this structure that makes it possible to move from a description of individual objects to a description of states and their quantitative occupations. In the proposed interpretation, this is the geometric expression of the inversion of quality and quantity: distinguishable modes define qualitative directions in state space, while occupation numbers determine the quantitative realization of the corresponding modes.

Thus, for Planck there are N resonators with mean energy U (continuous quality) and Pε energy quanta (discrete quantity), for which Planck wrote his principal equation $NU = P\varepsilon$. It is already clear from this equation that the more energetic quanta are fewer, and therefore their contribution to radiation is small, which makes it possible to avoid the theoretical ultraviolet catastrophe. The fact that for Planck the number P of quanta is primary, while N, the occupation number of this quantum, is secondary, is also reflected in the formula itself: it contains not (P−1), but specifically (N−1). So, $NU = P\varepsilon$ and it is beginning the combinations with repetition A.

In subsequent derivations of Bose–Einstein statistics, this emphasis is effectively shifted, and conversely, (P−1), rather than (N−1), appears. The numbers N and P are nevertheless very large, so this does not affect the calculations, but it is very important for understanding the essence of the matter.

**Conceptual and Mathematical Justification of the Method of Inversion of Qualities in Quantum Statistics**

The traditional physical paradigm relies on a corpuscular intuition in which the primary basis of reality is Quantity — a set of isolated objects (particles) localized in space. In this logic, the quantum state and the macroscopic structure of the system are regarded as a secondary Quality, emerging from the distribution of these objects over energy levels.

The method of inversion of qualities proposes a radical reversal of the familiar causal relation. Based on the dialectical idea of the primacy of Quality, one may postulate that the quantum state itself (the spatial-frequency mode of a field) is the initial physical reality. In such a paradigm, the quantum "building block" is defined not by a quantitative coordinate on a numerical axis, but by the modal question: "what quality is this?"

To describe the entire diversity of the microscopic world, the hypothetical coordinate system must become multidimensional, with each axis representing a fundamental, strictly discrete quality of matter:

- Frequency and energy axis (E=hω) — defines the dynamical quality of the quantum state.

- Mass axis (m) — fixes the quality of inertia or the “internal frequency” of the mode at rest.
- Charge axis (q) — represents the quality and symmetry of interaction (electromagnetic, color, etc.).
- Spin axis (s) — determines the fundamental spatial-quantum quality that separates matter into fermions and bosons.

An elementary particle does not bring these properties with it; rather, by occupying a definite point (the intersection of coordinates) in this space of qualities, the object is actualized and receives its physical “name.”

For example, a point combining the qualities of a definite frequency, zero mass, zero charge, and integer spin (s=1) gives the object the “name” photon. But if an object occupies a point with the quality of rest mass, negative charge, and half-integer spin (s=1/2), it is actualized under the “name” electron. At the same time, “names” can be scaled quantitatively (for example, $2h\omega_1$, with the other qualitative axes unchanged, means two photons of the same mode).

A physical confirmation of this logic is provided by black-body radiation. Here the general “name” photon conceals a whole continuous range along the frequency axis. A black body emits photons of many different frequencies, actualizing in space a moving spectrum of qualities differing in their energetic nature.

**1. Mathematical Apparatus: Operators as Tools for Counting Quality**

This philosophical inversion finds direct confirmation in the rigorous mathematical apparatus of quantum field theory — Dirac’s second-quantization method. In this scheme, the state of the system is specified not by particle coordinates, but by a vector in Fock space through the occupation numbers of modes. The principal mathematical tools here are creation and annihilation operators for a specific frequency mode.$a_\nu^\dagger a_\nu \nu$ .

The action of the creation operator on a quantum state with a spectral quality and a current occupation number is written as$\nu n_\nu$

$$a_\nu^\dagger|n_\nu\rangle = \sqrt{n+1}\,\lceil n_\nu + 1\rangle \quad .$$

Within the inversion approach, this formula acquires a fundamental meaning: the operator does not “bring” a material ball into space. It is a mathematical trigger that transfers the pre-existing Quality (the mode) to the next level of energetic activation.$a_\nu^\dagger \nu$

The particle-number (matter) operator itself is expressed as the product of the creation and annihilation operators:

$$\widehat{N}_\nu = a_\nu^\dagger \, a_\nu \ .$$

Its eigenvalue, when acting on the state vector, returns the occupation number:

$$\widehat{N}_\nu \lceil n_\nu \rangle = n_\nu \lceil n_\nu \rangle \ .$$

The formula clearly demonstrates that discrete matter in the form of particles is not the fundamental basis of reality, but an operational "counter" of energy, just as quantity is a counter of manifestations of a fundamental quality. A quantum particle is not a localized material object, but an operational marker, fixing the multiplicity of excitation of a specific mode.$n_\nu$

If a set of elements belongs to the same quantum mode, they are not a mechanical sum of independent objects, but represent a single Quality. A vivid physical embodiment of this principle is the phenomenon of the laser. A single photon is like a single sounding of the note "A" at a definite frequency. When the occupation number of this mode increases to millions, an ideal macroscopic unison — a laser beam — emerges. The laser demonstrates how a colossal Quantity of quantum objects is generated by the development and scaling of one and the same initial Quality.$n_{Ля} = 1( \, n_{Ля} \ \rightarrow 10^6 \, )$

**4. Resolution of Wave–Particle Duality through the Inversion of Qualities**

The century-old problem of wave–particle duality, which forced classical physicists to speak of the "split personality" of quantum objects (now a wave, now a particle), finds a natural and internally consistent logical resolution in the method of inversion of qualities.

In the traditional approach, the duality is perceived as a paradox: how can a localized, solid material "ball" (corpuscle) simultaneously possess the properties of an extended wave spread through space? Physicists were forced to introduce Born's probabilistic interpretation, asserting that the wave is merely an abstract "probability wave" for detecting the particle.

The method of inversion of qualities removes this paradox by eliminating the false dichotomy (separation) between wave and particle:

- Wave — pure Quality: An extended wave field in space (a mode with frequency ν) exists initially. It is a continuous, smooth structure of reality that determines the form and the rules of the game. The wave does not consist of particles; it is itself the fundamental qualitative state of matter.

- Particle — a step of Quantity (an act of counting): What we register with instruments as a "particle" (a detector click, a point on a photographic plate) is not the arrival of a solid microscopic object. It is the moment of quantized, discrete energy exchange between the wave mode and the instrument. It is a discrete "mark" in our counter, recording that the Quality has manifested itself exactly once. $n_\nu = 1$ )

  Thus, a quantum object does not "choose" whether to be a wave or a particle depending on the experiment. The object is always a wave in its ontological essence (Quality) and always manifests itself as a particle in the act of measurement (through Quantity).

In the double-slit experiment, the wave-Quality passes freely through both slits at once because the field mode is not localized. But when this wave encounters the screen, energy transfer occurs discretely — in portions hν. The screen acts as a matrix of "counters," recording individual manifestations of Quality. It seems to us that a point-particle has arrived, but this point is merely a local trace (a quantum imprint) of the extended spatial wave.

This view frees quantum mechanics from contrived mysticism and restores strict logic to it: the world is continuous in its Qualities (waves/modes) and discrete in its Quantities (particles/occupation numbers).

**2. Inversion of Qualities in Bose and Fermi Statistics: The Nature of the Counting Restriction**

The difference between the fundamental quantum statistics (Bose–Einstein and Fermi–Dirac) in the proposed concept is reduced to the internal capacity of Quality itself, that is, to the "design" of its counter.

- Bosons (symmetric states): For modes such as photons in a laser, the quality has unlimited quantitative capacity. The operators satisfy commutation relations. The occupation number may take any value: Quantity can infinitely scale a single Quality $[a_\nu, a_\mu^\dagger] = \delta_{\nu\mu} n_\nu = \{0, 1, 2, \dots . \infty\}$.
- Fermions (antisymmetric states): For modes such as electrons, Quality is alternative and modal. The operators obey anticommutation relations.$\{a_\nu, a_\mu^\dagger = \delta_{\nu\mu}$

  Because of this fundamental property, the creation operator applied twice to the same mode identically annihilates the state. This strict algebraic restriction gives rise to the famous Pauli exclusion principle. In the logic of inversion, the Pauli principle means that the fermionic "counter" of Quantity has only two divisions. A specific quality (an electron mode with a given set of quantum numbers) can either be absent from the physical world or be manifested strictly once. Fermionic matter structures the world precisely because its fundamental Qualities ("notes") prohibit quantitative duplication of themselves, so that the occupation number of each such state can take only the values $n_i=0$ or $n_i=1$.

**3. Cosmological Scale: Vacuum as the "Mode of All Modes"**

This approach is remarkably consonant with the concept of the distinguished physicist John Archibald Wheeler, "It from bit." Wheeler argued that all physical matter is fundamentally immaterial and arises from informational answers to binary questions. In the method of inversion, the occupation number of a mode functions precisely as Wheeler's information code: a specification of whether the quality is "on" and what the intensity of this activation is.

The logical completion of the inversion method is the concept of a unified physical vacuum as the "mode of all modes." The vacuum is not emptiness, but the universal ground state $|0\rangle$, an ocean of potential possibilities undergoing continuous zero-point fluctuations. Dirac's physical vacuum represents a total, undifferentiated proto-Quality containing the spectral potential of all wave modes of the Universe. The emergence of particular frequencies and the subsequent growth of their occupation numbers is the fractal unfolding of the primordial Quality. Physical reality, in this reading, is not a mosaic of material building blocks, but a polyphony of coexisting qualitative states of a unified substance.

1.

**Part 1. Two-Level System: Birth of Quantum Statistics from the Inversion of Qualities**

Let us apply the method of inversion to the simplest physical model: a system having only two possible states — "this level" (the selected state) and "not this level" (the rest of the state space).

In the classical representation, we would take an ensemble of N material objects and begin distributing them among these two levels, counting combinations. This is exactly what Boltzmann does, calling combinations without repetition "complexions." The method of inversion requires us to reverse the problem: we fix the quantum state itself as the fundamental Quality, while physical objects are reinterpreted as Quantity (the occupation number) of this quality. Since quantum particles are absolutely indistinguishable, they have no "personal numbers." We cannot ask which particular object occupied the level. We can only determine in what quantity a given quality is manifested in the system. The mathematical analysis of this inversion leads to two fundamental scenarios — the two quantum statistics.

**1. Fermions: Combinations without Repetition**

The first scenario arises when nature imposes a strict restriction on the manifestation of quality: it either exists or it does not. In such a system, the occupation number (“number of objects”) for a particular level can take only two values: 0 (quality absent) or 1 (quality manifested).

This is exactly what Fermi does. His statistics is constructed in this way: the possible occupations of quantum levels are considered. Now for Fermi, N is the number of levels, and these levels, already as occupation numbers, represent Quantity, while the object — whether it is present at this level or absent — becomes the quality of the level. Fermi points his “finger” not at the object, as Boltzmann does, asking whether it is selected or not, but at the level, asking whether it is occupied or not. In combinatorics, this choice is described as combinations without repetition. We combine not objects, but the quantum qualities themselves, and none of them can be selected twice. This strict prohibition on repetition of qualities physically represents the Pauli exclusion principle [11]. If we try to combine many such systems, they must occupy increasingly higher and new energy levels, forming a complex and ordered macroscopic structure of matter. From this fermionic logic of “combinations without repetition” arise chemistry as a whole, the structure of atomic electron shells, and the periodic table of Mendeleev. Thus, in Fermi–Dirac statistics, the occupation number (0 or 1) is transformed into a qualitative characteristic of the level — its “activity” or “emptiness.” The number of levels N, meanwhile, acts as the quantity of available qualities, while the occupation number itself functions as a binary switch determining whether the given quality is manifested in reality. But there is another statistic, the statistics of bosons. This is the one Planck used for black-body radiation.

### 2. Bosons: Combinations with Repetition

The second scenario removes the restriction on uniqueness. The same quantum quality may be selected by the system repeatedly. The occupation number (“number of objects”) at the selected level becomes arbitrary: $n = 0, 1, 2, 3\ldots.\infty$

In combinatorial terms, these are pure combinations with repetition. The system can accumulate quantity in a cascade within the same unchanged quality.

In the physical world, such particles are called bosons. The ability of bosons to “repeat” one another’s quality leads to their tendency to collect collectively in the same quantum state. At extreme cooling this produces Bose–Einstein condensation, while in the electromagnetic field it produces photon synergy, making a laser beam possible.

## Part 2. Harmonic Oscillator and Derivation of Planck’s Formula: The Magic of Equidistant Systems

Let us transfer the inversion method to an infinite system — the quantum harmonic oscillator. This model serves as an ideal bridge between the abstract “qualities” of fields and the “quantity” of particles because its energy spectrum possesses a unique property — equidistance (equal spacing of levels).

The energy of each oscillator level is given by the formula:

$$E_n = h\omega\left(n + \frac{1}{2}\right) \ , \qquad \text{где } n = 1,2,3 \ldots \infty \ .$$

If the zero-point energy is omitted, the distance between any two neighboring levels is strictly the same and equal to the fundamental quantum of energy.$\frac{1}{2}h\omega\Delta E = h\omega$ .

**Step 1. Logical Inversion: The Level Number as the Number of Objects**

It is precisely equidistance that makes it possible to perform a flawless logical inversion. In the classical paradigm we would say: “The oscillator has moved to the n-th energy level.”

The method of inversion of qualities proposes looking at this differently: we fix the wave mode of frequency as a single Quality. Then the transition of the system to the n-th level means that this Quality has generated within itself a Quantity equal to n identical, absolutely indistinguishable quanta (photons) with the corresponding energy.$\omega h\omega$

The level number of the harmonic oscillator is physically identical to the occupation number (the number of objects) in this mode. Since photons are free to be created in any volume and are not restricted by the Pauli principle, we are dealing with combinations with repetition.

**Step 2. Derivation of Planck’s Denominator**

To demonstrate the operational power of the method of inversion of qualities, let us show how the standard mathematical apparatus of the harmonic oscillator can be reinterpreted through the combinatorics of modes rather than the distribution of particles.
Under thermodynamic equilibrium at temperature T, the probability of finding the system in a state with energy proportional to the classical Gibbs factor is: $E_n = h\omega P_n \approx e^{-\frac{nh\omega}{k_B T}}$

Let us find the mean number of quanta (the mean oscillator level number) as the mathematical expectation over all possible occupation configurations from 0 to $\langle n \rangle \infty$ .

$$\langle n \rangle = \frac{\sum_0^\infty n e^{-\frac{nh\omega}{k_B T}}}{\sum_0^\infty e^{-\frac{nh\omega}{k_B T}}} \ .$$

For convenience, introduce the notation. Then the expression becomes a ratio of two series:

$$x = e^{-\frac{h\omega}{k_B T}} .$$

The denominator is an infinite decreasing geometric progression:

$$\sum_{n=0}^{\infty} x^n = \frac{1}{1-x} .$$

The numerator is the derivative of this series multiplied by x:

$$\sum_{n=0}^{\infty} n x^n = x \frac{d}{dx}\left(\frac{1}{1-x}\right) = \frac{x}{(1-x)^2} .$$

Substituting these sums back into the formula for $\langle n \rangle$ :

$$\langle n \rangle = \frac{\frac{x}{(1-x)^2}}{\frac{1}{1-x}} = \frac{x}{1-x} = \frac{1}{\frac{1}{x} - 1} .$$

Returning to the substitution, we obtain the fundamental bosonic Planck denominator:

$$x = e^{-\frac{h\omega}{k_B T}}$$

$$\langle n \rangle = \frac{1}{e^{\frac{h\omega}{k_B T}} - 1} .$$

From the standpoint of the inversion of qualities, this classical result (1) acquires a new meaning. Planck's denominator is not a law describing the distribution of a "crowd" of distinguishable particles. It is a strict mathematical measure of how a single spatial mode (Quality) unfolds its internal quantitative potential under the influence of thermal chaos.

This derivation of Planck's formula can also be obtained by agreeing with Planck that between neighboring states of one frequency mode there is the same energy quantum, and by noting that the probability of occupation of the mode is proportional to the relative number of quanta in this mode.

**Step 3. Escape from the "Ultraviolet Catastrophe"**

Multiplying the mean number of quanta by the energy of one quantum and by the geometrical density of electromagnetic modes in the cavity, we obtain the complete Planck formula for black-body radiation: $\langle n \rangle h\omega(\frac{\omega^3}{\pi^2 c^3})$

$$u(\omega,T) = \frac{\omega^3}{\pi^2 c^3} \frac{1}{e^{-\frac{nh\omega}{k_B T}} - 1} \ .$$

Classical physics (the Rayleigh–Jeans law) required each mode to carry the mean energy of thermal chaos. But this led to infinite radiation energy — the "ultraviolet catastrophe" [4, 7, 9].$u(\omega,T) = \frac{\omega^3}{\pi^2 c^3}\ k_B T\ )k_B T\omega \to \infty$

Our inversion technique clearly shows how nature solves this problem combinatorially. At high frequencies, the energy of one quantum becomes enormous compared with the thermal energy. The exponential in the denominator grows rapidly, and the probability of creating even one quantum (realizing such a "quantity" within the given "quality") falls to zero. The bosonic minus term "switches off" high-frequency oscillations, saving the spectrum from infinity and yielding a finite integral relation proportional to $T^4$ (the Stefan-Boltzmann law).

**Part 3. Extrapolation of the Method: Limits of Applicability and Geometry of State Spaces**

The method of inversion of qualities clearly demonstrates that the macroscopic behavior of a system is predetermined by the internal geometry of its energy levels. Consideration of two contrasting examples — a crystal lattice and the hydrogen atom — makes it possible to draw a clear boundary of applicability for this logic.

**1. Thermal Vibrations of Crystals (Einstein Model)**

In 1907 Albert Einstein applied an approach fully equivalent to our method to the thermodynamics of solids [6]. Before him, classical physics asserted that the heat capacity of metals should remain constant at all temperatures (the Dulong–Petit law), which was in serious contradiction with experiments under strong cooling.

Einstein proposed treating the crystal as an ensemble of atoms oscillating independently at a fixed frequency.$\omega$

- Quality: Each axis of atomic vibration in the lattice is fixed as a quantum Quality (a mode of the sound field).

- Quantity: The spacing between the energy levels of such an oscillation is strictly equidistant. Consequently, excitation of the oscillation is the creation of faceless, indistinguishable quanta of sound — phonons.

  Since phonons have no restriction on the number of quanta occupying the cells, we again have combinations with repetition. The inversion method works here without modification: the mean number of sound quanta in a mode is given by the same bosonic Planck denominator.

$$\langle n \rangle = \frac{1}{e^{\frac{nh\omega}{k_B T}} - 1} \quad .$$

As the temperature decreases, the thermal energy becomes smaller than the sound quantum. Phonons cease to be created, atomic vibrations “freeze,” and the heat capacity of the crystal rapidly falls to zero. Our inversion model transfers from light to sound without any change, confirming its universality for all equidistant systems.$T \to 0h\omega$

**2. Structure of the Hydrogen Atom: The Nonequidistant Limit**

The logic that “the level number equals the number of objects within the quality” changes radically when we move to systems with a curved, nonuniform geometry of levels. The textbook example is the hydrogen atom.

The energy of its levels is given by Bohr’s formula:

$$E_n = -\frac{13{,}6\ \text{эВ}}{n^2}\ , \qquad \text{где}\ \ n = 1{,}2{,}3\ldots\infty\,.$$

Here the distances between neighboring levels differ: the step from the first level to the second is enormous (10.2 eV), while the steps between upper levels become increasingly smaller.

1. Breakdown of Equivalence: The level number n is no longer identical to the occupation number. When the atom is on the third level, it does not contain “three identical quanta.” The energy is not a multiple of the quantum energy.$E_3 E_1$
2. Transition to Individual States: The electron inside the atom is a single fermion. It cannot multiply or accumulate within one level.
3. In a thermodynamic ensemble of such atoms, the “qualities” become the levels themselves, but for each particular atom this is a strictly nonrepeating choice (one atom — one level).$E_n$

   Therefore, the mean population of levels in an ensemble of hydrogen atoms obeys not Planck’s denominator but the classical Maxwell–Boltzmann distribution. Without equidistant levels, the bosonic magic disappears, and the system returns to counting individual, fixed states. $\langle n \rangle \approx e^{-\frac{nh\omega}{k_B T}}$

**Conclusion: The Architecture of Reality in the Language of Occupation Numbers**

The method of inversion of qualities proposed and developed in this work takes us beyond a simple combinatorial trick. It exposes the very logic by which quantum field theory is structured. Modern particle physics has finally abandoned the concept of “material balls” moving through empty space. At the fundamental level of the Universe there are no individual objects — there are only unified quantum fields (electron, photon, Higgs), permeating space. Each such field is a global Quality. What we register in the macroscopic world as “matter” or “radiation” is merely local Quantity, a degree of excitation (occupation number) of this field.

By reversing Hegel’s familiar triad and placing Quality before Quantity, we have discovered a common root of quantum symmetry:

- If the Universe imposes a prohibition on the repetition of qualities, we obtain the world of fermions (combinations without repetition). This rigid combinatorial framework prevents electrons from collapsing into the nucleus, arranges the elements in the periodic table, and creates the stable solid matter of which we ourselves are made.
- If the Universe permits infinite repetition of qualities, we obtain the world of bosons (combinations with repetition). In equidistant systems (harmonic oscillators), this choice gives rise to photons and phonons. It allows energy to accumulate freely at one quantum point, igniting stars, generating laser beams, and gathering billions of atoms into a single macroscopic Bose-condensed state [5–8].

Thus, all the complexity and diversity of the surrounding world — from solid silicon to sunlight — is a direct consequence of how nature combines its fundamental qualities.